\documentclass[a4paper,12pt,oneside]{article}
\usepackage{amsmath,amsfonts,amssymb,amsmath,latexsym,amsthm,enumerate}

\newtheorem{thm}{Theorem}
\newtheorem{cor}[thm]{Corollary}
\theoremstyle{definition}

\theoremstyle{plain}

\begin{document}
\title {Some identities of higher-order Euler polynomials arising from Euler basis}
\author{by \\Dae San Kim and Taekyun Kim}\date{}\maketitle

\begin{abstract}
\noindent The purpose of this paper is to present a systematic study of some families of higher-order Euler numbers and polynomials. In particular, by using the basis property of higher-order Euler polynomials for the space of polynomials of degree less than and equal to $n$, we derive some interesting identities for the higher-order Euler polynomials.
\end{abstract}

\section  {Introduction}
As is well known, the $n$-th Euler polynomials of order $r$ are defined by the generating function to be
\begin{equation}\label{inteqn1}(\frac{2}{e^{t}+1})^{r}e^{xt}=e^{E^{(r)}(x)t}=\sum_{n=0}^{\infty}E_{n}^{(r)}(x)\frac{t^{n}}{n!}\quad(r\in\mathbf{Z}_{+})\,,\end{equation}with the usual convention about replacing $(E^{(r)}(x))^{n}$ by $E_{n}^{(r)}(x)$\quad(see\,\,\lbrack1-11\rbrack)\,.\\
In the special case, $x=0,\, E_{n}^{(r)}(0)=E_{n}^{(r)}$ are called the $n$-th Euler numbers of order $r$.\\
By(\ref{inteqn1}),\quad we easily get
\begin{align}\label{inteqn2} E_{n}^{(r)}(x) &= \sum_{l=0}^{n}\left(\begin{array}{c}n\\l\end{array}\right)E_{l}^{(r)}x^{n-l}=\sum_{l=0}^{n}\left(\begin{array}{c}n\\l\end{array}\right)E_{n-l}^{(r)}x^{l}\\
                                           &=\sum_{n=n_{1}+\cdots+n_{r}+n_{r+1}}\left(\begin{array}{c}n\\n_{1},\cdots,n_{r},n_{r+1}\end{array}\right)E_{n_{1}}E_{n_{2}}\cdots E_{n_{r}}x^{n_{r+1}}\,.\nonumber
\end{align}
From(\ref{inteqn2}),\quad we note that the leading coefficient of $E_{n}^{(r)}(x)$ is given by
\begin{equation}\label{inteqn3} \sum_{n_{1}+\cdots+n_{r}=0}\left(\begin{array}{c}n\\n_{1},\cdots,n_{r}\end{array}\right)E_{n_{1}}E_{n_{2}}\cdots E_{n_{r}}=1\,.
\end{equation}
Thus,\, $E_{n}^{(r)}(x)$ is a monic polynomial of degree $n$ with rational coefficients.\\
From(\ref{inteqn1}),\, we have $E_{n}^{(0)}(x)=x^{n}$.\,\,It is not difficult to show that
\begin{eqnarray}\label{inteqn4}
\frac{dE_{n}^{(r)}(x)}{dx}=nE_{n-1}^{(r)}(x)\,,\quad E_{n}^{(r)}(x+1)+E_{n}^{(r)}(x)=2E_{n}^{(r-1)}(x), \quad (\text{see\,\lbrack11-18\rbrack})\,.
\end{eqnarray}
Now, we define two linear operators $\tilde{\triangle}$ and $D$ on the space of real-valued differentiable functions as follows:
\begin{eqnarray}\label{inteqn5}
\tilde{\triangle}f(x)=f(x+1)+f(x),\quad Df(x)=\frac{df(x)}{dx}\,.
\end{eqnarray}
Then we see that $\tilde{\triangle}D=D\tilde{\triangle}$\,.\\
\\
Let $V_{n}=\lbrace p\,(x)\in\mathbf{Q}[x]\vert$ deg $p\,(x)\leq n\rbrace$ be the ($n+1$)-dimensional vector space over $\mathbf{Q}$.\,\,Probably, $\lbrace 1, x, \cdots, x^{n}\rbrace$ is the most natural basis for $V_{n}$.\, But $\lbrace E_{0}^{(r)},E_{1}^{(r)},\cdots,E_{n}^{(r)}\rbrace$ is also a good basis for the space $V_{n}$ for our purpose of arithmetical and combinatorial applications of the higher-order Euler polynomials.\\
If $p\,(x)\in V_{n}$, then $p(x)$ can be expressed by
\begin{equation*}p\,(x)=b_{0}E_{0}^{(r)}(x)+b_{1}E_{1}^{(r)}(x)+\cdots+b_{n}E_{n}^{(r)}(x)\,.\end{equation*}
In this paper, we develop methods for computing $b_{l}$ from the information of $p\,(x)$ and apply those results to arithmetically and combinatorially interesting identities involving $E_{0}^{(r)},E_{1}^{(r)},\cdots,E_{n}^{(r)}$\,.
\section{Higher-order Euler polynomials}
From(\ref{inteqn5}),\quad we have
\begin{equation}\label{inteqn6}\tilde{\triangle}E_{n}^{(r)}(x)=E_{n}^{(r)}(x+1)+E_{n}^{(r)}(x)=2E_{n}^{(r-1)}(x)\,,\end{equation}
and
\begin{equation}\label{inteqn7}DE_{n}^{(r)}(x)=nE_{n-1}^{(r)}(x)\,.\end{equation}
Let us assume that $p\,(x)\in V_{n}$.\, Then $p\,(x)$ can be generated by $ E_{0}^{(r)}(x),E_{1}^{(r)}(x),\\\cdots,E_{n}^{(r)}(x)$ to be
\begin{equation}\label{inteqn8}p\,(x) = \sum_{k=0}^{n}b_{k}E_{k}^{(r)}(x)\,.\end{equation}
Thus, by({\ref{inteqn8}), we get
\begin{equation*}\tilde{\triangle}p\,(x)=\sum_{k=0}^{n}b_{k}\tilde{\triangle}E_{k}^{(r)}(x)=2\sum_{k=0}^{n}b_{k}E_{k}^{(r-1)}(x)\,,\end{equation*}
and
\begin{equation*}\tilde{\triangle}^{2}p\,(x)=2\sum_{k=0}^{n}b_{k}\tilde{\triangle}E_{k}^{(r-1)}(x)=2^{2}\sum_{k=0}^{n}b_{k}E_{k}^{(r-2)}(x)\,.\end{equation*}
Continuing this process, we have
\begin{equation}\label{inteqn9}\tilde{\triangle}^{r}p(x)=2^{r}\sum_{k=0}^{n}b_{k}E_{k}^{(0)}(x)=2^{r}\sum_{k=0}^{n}b_{k}x^{k}\,.\end{equation}
Let us take the operator $D^{k}$ on (\ref{inteqn9}).\,\,Then
\begin{align}\label{inteqn10}D^{k}\tilde{\triangle}^{r}p(x)&=2^{r}\sum_{l=k}^{n}b_{l}l(l-1)\cdots(l-k+1)x^{l-k}\\
                                                           &=2^{r}\sum_{l=k}^{n}b_{l}\frac{l!}{(l-k)!}x^{l-k}\nonumber\\
                                                           &=2^{r}\sum_{l=k}^{n}b_{l}k!\left(\begin{array}{c}l\\k\end{array}\right)x^{l-k}\nonumber.
\end{align}
Let us  take $x=0$ on (\ref{inteqn10}).\,\, Then we get
\begin{equation}\label{inteqn11}D^{k}\tilde{\triangle}^{r}p\,(0)=2^{r}b_{k}k!\,.\end{equation}
From (\ref{inteqn11}),\quad we have
\begin{align}\label{inteqn12} b_{k}&=\frac{1}{2^{r}k!}D^{k}\tilde{\triangle}^{r}p\,(0)=\frac{1}{2^{r}k!}\tilde{\triangle}^{r}D^{k}p\,(0)\\
     &=\frac{1}{2^{r}k!}\sum_{j=0}^{r}\left(\begin{array}{c}r\\j\end{array}\right)D^{k}p\,(j)\,.\nonumber
\end{align}
Therefore, by (\ref{inteqn8})and(\ref{inteqn12}), we obtain the following theorem.
\begin{thm}\label{Thm1}
For $n, r\in\mathbf{Z}_{+},\quad p(x)\in V_{n}$, we have
\begin{equation*}p\,(x)=\frac{1}{2^{r}}\sum_{k=0}^{n}(\sum_{j=0}^{r}\frac{1}{k!}\left(\begin{array}{c}r\\j\end{array}\right)D^{k}p\,(j))E_{k}^{(r)}(x)\,.\end{equation*}
\end{thm}
Let us take $p\,(x)=x^{n}\in {V}_{n}$.\, Then we easily see that $D^{k}x^{n}=\frac{n!}{(n-k)!}x^{n-k}$\,.\\
Thus, by Theorem \ref{Thm1}, we get
\begin{align}\label{inteqn13}x^{n}&=\frac{1}{2^{r}}\sum_{k=0}^{n}\sum_{j=0}^{r}\frac{1}{k!}\left(\begin{array}{c}r\\j\end{array}\right)\frac{n!}{(n-k)!}j^{n-k}E_{k}^{(r)}(x)\\
                                     &=\frac{1}{2^{r}}\sum_{k=0}^{n}\sum_{j=0}^{r}\left(\begin{array}{c}r\\j\end{array}\right)\left(\begin{array}{c}n\\k\end{array}\right)j^{n-k}E_{k}^{(r)}(x)\,.\nonumber
\end{align}
Therefore, by (\ref{inteqn13}), we obtain the following corollary.
\begin{cor} For $n, r\in\mathbf{Z}_{+}$, we have
\begin{equation*}x^{n}=\frac{1}{2^{r}}\sum_{k=0}^{n}\sum_{j=0}^{r}\left(\begin{array}{c}r\\j\end{array}\right)\left(\begin{array}{c}n\\k\end{array}\right)j^{n-k}E_{k}^{(r)}(x)\,.\end{equation*}
\end{cor}
Let $p\,(x)=B_{n}^{(s)}(x) (s\in\mathbf{Z}_{+})$.\,\,Then we have
\begin{equation}\label{inteqn14}D^{k}B_{n}^{(s)}(x)=\frac{n!}{(n-k)!}B_{n-k}^{(s)}(x)\,.\end{equation}
By Theorem \ref{Thm1}, we get
\begin{equation}\label{inteqn15}B_{n}^{(s)}(x)=\frac{1}{2^{r}}\sum_{k=0}^{n}\sum_{j=0}^{r}\left(\begin{array}{c}r\\j\end{array}\right)\left(\begin{array}{c}n\\k\end{array}\right)B_{n-k}^{(s)}(j)E_{k}^{(r)}(x)\,.\end{equation}
Therefore, by(\ref{inteqn15}), we obtain the following corollary.
\begin{cor} For $n, s, r\in\mathbf{Z}_{+}$,\quad we have
\begin{equation*}B_{n}^{(s)}(x)=\frac{1}{2^{r}}\sum_{k=0}^{n}\sum_{j=0}^{r}\left(\begin{array}{c}r\\j\end{array}\right)\left(\begin{array}{c}n\\k\end{array}\right)B_{n-k}^{(s)}(j)E_{k}^{(r)}(x)\,,\end{equation*}
where $B_{n}^{(s)}(x)$ are the $n$-th Bernoulli polynomials of order $s$.
\end{cor}
\noindent It is well known that
\begin{eqnarray}\label{inteqn16}\frac{t}{e^{t}-1}e^{xt}=\sum_{n=0}^{\infty}B_{n}(x)\frac{t^{n}}{n!}\,,\quad \frac{2}{e^{t}+1}e^{xt}=\sum_{n=0}^{\infty}E_{n}(x)\frac{x^{n}}{n!}\,.
\end{eqnarray}
In the special case, $x=0$, let $B_{n}(0)=B_{n}$,\,\,$E_{n}(0)=E_{n}$\,.
From(\ref{inteqn16}), we easily derive the following identity:
\begin{equation}\label{inteqn17}B_{n}(x)=\sum_{k=0,k\neq 1}^{n}\left(\begin{array}{c}n\\k\end{array}\right) B_{k}E_{n-k}(x)\in V_{n}\,.\end{equation}
Let us take $p\,(x)=B_{n}(x)$.\,\,Then we have
\begin{equation}\label{inteqn18}D^{k}B_{n}(x)=n(n-1)\cdots(n-k+1)B_{n-k}(x)=\frac{n!}{(n-k)!}B_{n-k}(x)\,.\end{equation}
Therefore, by Theorem \ref{Thm1},\,\,(\ref{inteqn17}) and (\ref{inteqn18}), we obtain the following theorem.
\begin{thm}
For $n, r\in\mathbf{Z}_{+}$,\quad we have
\begin{equation*}\sum_{k=0,k\neq1}^{n}\left(\begin{array}{c}n\\k\end{array}\right)B_{k}E_{n-k}(x)=\frac{1}{2^{r}}\sum_{k=0}^{n}\sum_{j=0}^{r}\left(\begin{array}{c}r\\j\end{array}\right)\left(\begin{array}{c}n\\k\end{array}\right)B_{n-k}(j)E_{k}^{(r)}(x)\,.\end{equation*}
\end{thm}
\quad
\\
\\
Let us consider $p\,(x)=\sum_{k=0}^{n}B_{k}(x)B_{n-k}(x)$\,.\\
Then we have
\begin{equation}\label{inteqn19}
D^{k}p\,(x)=\frac{(n+1)!}{(n-k+1)!}\sum_{l=k}^{n}B_{l-k}(x)B_{n-l}(x)\,.
\end{equation}
Thus, by Theorem \ref{Thm1} and (\ref{inteqn19}), we obtain the following theorem.\\
\\
\begin{thm}
For $r, n\in\mathbf{Z}_{+}$, we have
\begin{equation*}\sum_{k=0}^{n}B_{k}(x)B_{n-k}(x)=\frac{1}{2^{r}}\sum_{k=0}^{n}\sum_{l=k}^{n}\sum_{j=0}^{r}\left(\begin{array}{c}r\\j\end{array}\right)\left(\begin{array}{c}n+1\\k\end{array}\right)B_{l-k}(j)B_{n-l}(j)E_{k}^{(r)}(x)\,.\end{equation*}
\end{thm}
\quad
\\
\\
Let $n,\,m\in\mathbf{Z}_{+}$\,, with $n\geq m+2$.\quad Then we have
\\
\begin{align}\label{inteqn20}
&B_{m}(x)B_{n-m}(x)\nonumber\\
&=\sum_{l=0}^{\infty}\lbrace\left(\begin{array}{c}m\\2l\end{array}\right)(n-m)+\left(\begin{array}{c}n-m\\2l\end{array}\right)m\rbrace\frac{B_{2l}B_{n-2l}(x)}{n-2l}+(-1)^{m+1}\frac{B_{n}}{\left(\begin{array}{c}n\\m\end{array}\right)}\,.\end{align}
\\
Let us take $p\,(x)=B_{m}(x)B_{n-m}(x)\in V_{n}$\,.\\
Then we have
\begin{equation*}\tag{21}
\begin{split}
&D^{k}p\,(x)\\
&=\sum_{l=k}^{\infty}\biggr{\{}\binom{m}{2l}(n-m)+\binom{n-m}{2l}m\biggr{\}}\frac{B_{2l}}{n\!\!-\!\!2l}\times\frac{(n-2l)!}{(n-2l-k)!}B_{n-2l-k}(x).
\end{split}
\end{equation*}
Therefore, by Theorem \ref{Thm1} and (21), we obtain the following theorem.
\begin{thm}\label{Thm6}
For $n,m\in\mathbf{Z}_{+}$ with $n\geq m+2$,\quad we have
\begin{equation*}
\begin{split}
&B_{m}(x)B_{n-m}(x)\\
&=\frac{1}{2^{r}}\sum_{k=0}^{n}\Biggr{\{}\sum_{l=k}^{\infty}\sum_{j=0}^{r}\binom{r}{j}\binom{n-2l}{k}
\\&\times\biggr{(}\binom{m}{2l}(n-m)+\binom{n-m}{2l}m\biggr{)}
\frac{B_{2l}B_{n-2l-k}(j)}{n-2l}\Biggr{\}}E_{k}^{(r)}(x).
\end{split}
\end{equation*}
 \end{thm}
\quad
\\
\\
${\large\mathbf{Remark}}$.
By using Theorem \ref{Thm1}, we can find many interesting identities related to Bernoulli and Euler polynomials.\nonumber

\quad
\\
\\
\\

\author{Department of Mathematics, Sogang University, Seoul 121-742, Republic of Korea
\\e-mail: dskim@sogang.ac.kr}\\
\\
\author{Department of Mathematics, Kwangwoon University, Seoul 139-701, Republic of Korea
\\e-mail: tkkim@kw.ac.kr}
\end{document}